\documentclass[12pt]{article}
\usepackage{amssymb}
\usepackage{latexsym,bm}
\usepackage{graphicx}
\usepackage{amsmath}
\usepackage{mathrsfs}

\newcommand{\bea}{\begin{eqnarray*}}
\newcommand{\eea}{\end{eqnarray*}}
\newcommand{\be}{\begin{equation}}
\newcommand{\ee}{\end{equation}}
\newcommand{\ben}{\begin{eqnarray*}}
\newcommand{\een}{\end{eqnarray*}}

\date{}
\begin{document}
\title{The Tur\'{a}n number of book graphs\footnote{E-mail addresses:
{\tt mathyjr@163.com}(J.Yan),
{\tt zhan@math.ecnu.edu.cn}(X.Zhan).}}
\author{\hskip -10mm Jingru Yan and Xingzhi Zhan\thanks{Corresponding author.}\\
{\hskip -10mm \small Department of Mathematics, East China Normal University, Shanghai 200241, China}}\maketitle
\begin{abstract}
 Given a graph $H$ and a positive integer $n,$ the Tur\'{a}n number of $H$ for the order $n,$ denoted
 ${\rm ex}(n,H),$ is the maximum size of a simple graph of order $n$ not containing $H$ as a subgraph. The book with $p$ pages, denoted $B_p$, is the graph that consists of $p$ triangles sharing a common edge.  Bollob\'{a}s and Erd\H{o}s initiated the research on the Tur\'{a}n number of book graphs in 1975.
 The two numbers ${\rm ex}(p+2,B_p)$ and  ${\rm ex}(p+3,B_p)$ have been determined by Qiao and Zhan.
 In this paper we determine the numbers ${\rm ex}(p+4,B_p),$ ${\rm ex}(p+5,B_p)$ and ${\rm ex}(p+6,B_p),$ and characterize the corresponding extremal graphs
 for the numbers ${\rm ex}(n,B_p)$ with $n=p+2,\,p+3,\,p+4,\,p+5.$
\end{abstract}

{\bf Key words.} Tur\'{a}n number; book; triangle; extremal graph

{\bf Mathematics Subject Classification.} 05C35, 05C75

\section{Introduction}

We consider finite simple graphs. The {\it order} of a graph $G,$ denoted $|G|,$ is its number of vertices, and the {\it size} its number of edges.

{\bf Definition 1.} Given a graph $H$ and a positive integer $n,$ the {\it Tur\'{a}n number of $H$ for the order $n,$} denoted ${\rm ex}(n,H),$
is the maximum size of a simple graph of order $n$ not containing $H$ as a subgraph.

Determining the Tur\'{a}n number of various graphs $H$ is one of the main topics in extremal graph theory [1].

{\bf Definition 2.} The {\it book} with $p$ pages, denoted $B_p$, is the graph that consists of $p$ triangles sharing a common edge.

$B_5$ is depicted in Figure 1. Note that the graph $B_p$ has order $p+2.$
\vskip 3mm
\par
 \centerline{\includegraphics[width=2in]{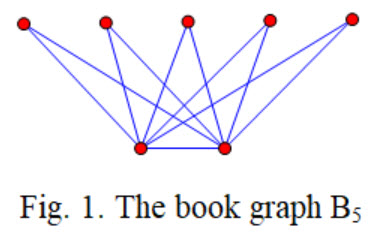}}
\par

We denote by $K_p$ the complete graph of order $p.$  $B_1$ is just $K_3.$
Denote by $K_{s,t}$ the complete bipartite graph whose partite sets have cardinalities $s$ and $t.$
A classic result of Mantel [5] from 1907 states that  ${\rm ex}(n, K_3)=\lfloor n^2/4\rfloor$ and the balanced complete
bipartite graph $K_{\lfloor n/2\rfloor,\lceil n/2\rceil}$ is the unique extremal graph.

In 1975 Bollob\'{a}s and Erd\H{o}s [2] posed the conjecture that ${\rm ex}(n,B_{\lceil n/6\rceil})\le n^2/4,$ which was proved
 by Edwards [3] and independently by Had\v{z}iivanov and Nikiforov [4].

The two numbers ${\rm ex}(p+2,B_p)$ and  ${\rm ex}(p+3,B_p)$ have been determined by Qiao and Zhan [6] in solving a problem of Erd\H{o}s
and they posed the problem [6, Problem 3] of determining ${\rm ex}(n,B_p)$ for a general order $n.$
In this paper we determine the numbers ${\rm ex}(p+4,B_p),$ ${\rm ex}(p+5,B_p)$ and ${\rm ex}(p+6,B_p)$ and characterize the corresponding extremal graphs
for the numbers ${\rm ex}(n,B_p)$ with $n=p+2,\,p+3,\,p+4,\,p+5.$

For graphs we will use equality up to isomorphism, so $G=H$ means that $G$ and $H$ are isomorphic.
Given graphs $G$ and $H,$ the notation $G+H$ means the disjoint union of $G$ and $H,$ and $G\vee H$ denotes the {\it join} of $G$ and $H,$ which is
obtained from $G+H$ by adding edges joining every vertex of $G$ to every vertex of $H.$  $m H$ denotes the disjoint union of $m$ copies of a graph $H.$
We denote by $V(G)$ and $E(G)$ the vertex set and edge set of a graph $G$ respectively. $\overline{G},$ $\delta (G)$ and $\Delta (G)$ will denote the
complement, minimum degree and maximum degree of a graph $G$ respectively. For $A\subseteq V(G),$ we denote by $G[A]$ the subgraph of $G$ induced by $A.$

For a vertex $v\in V(G),$ $N_G(v)$ and ${\rm deg}_G(v)$ will denote the neighborhood and degree of $v$ in $G$ respectively. If the graph $G$
 is clear from the context, we will omit it as the subscript. A vertex of degree $0$ is called an {\it isolated vertex} and a vertex of degree
 $1$ a {\it leaf.}

 We denote by $C_s$ and $P_t$ the cycle of length $s$ and path of order $t$ respectively.

{\bf Notation.} For an even positive integer $n,$ the notation $K_n-PM$ denotes the graph obtained from the complete graph $K_n$
by deleting all the edges in a perfect matching of $K_n$; i.e., it is the complement of $(n/2)K_2.$

\section{Main Results}

We will repeatedly use the first two lemmas below.

{\bf Lemma 1.} {\it Let $n$ and $p$ be positive integers with $n\ge p+2.$ Then a graph $G$ of order $n$ does not contain $B_p$ as a subgraph
if and only if for every edge $xy\in E(G),$ $|N_{\overline{G}}(x)\cup N_{\overline{G}}(y)|\ge n-p-1.$ Consequently, if $G$ is such a graph then
$\overline{G}$ has at most one isolated vertex.}

{\bf Proof.} For $xy\in E(G),$ $V(G)\setminus \{x,y\}=(N_G(x)\cap N_G(y))\cup (N_{\overline{G}}(x)\cup N_{\overline{G}}(y))$
and $(N_G(x)\cap N_G(y))\cap (N_{\overline{G}}(x)\cup N_{\overline{G}}(y))=\phi.$ Hence
$$
|N_{\overline{G}}(x)\cup N_{\overline{G}}(y)|=n-2-|N_G(x)\cap N_G(y)|.
$$
$G$ does not contain $B_p$ if and only if for every edge $xy\in E(G),$ $|N_G(x)\cap N_G(y)|\le p-1$ if and only if
$|N_{\overline{G}}(x)\cup N_{\overline{G}}(y)|\ge n-p-1.$ This shows the first conclusion. The second conclusion is obvious and it
also follows from the first immediately. $\Box$

{\bf Lemma 2.} {\it Let $n$ and $p$ be positive integers with $n\ge p+5$ and let $G$ be a graph of order $n.$ If $\overline{G}$ has a component
which is a cycle of length at least $4,$ then $G$ contains $B_p.$ }

{\bf Proof.} Let $C$ be a component of $\overline{G}$ which is a cycle of length at least $4,$ and let $x,\,z,\,y$ be three consecutive vertices of
$C.$ Then $xy\in E(G)$ and $|N_{\overline{G}}(x)\cup N_{\overline{G}}(y)|\le 3<n-p-1,$ since $n\ge p+5.$ By Lemma 1, $G$ contains $B_p.$ $\Box$

{\bf Lemma 3.} (Qiao and Zhan [6]) {\it If $p$ is an even positive integer, then ${\rm ex}(p+2,B_p)=p(p+2)/2$ and ${\rm ex}(p+3,B_p)=p(p+4)/2;$ if $p$ is an
odd positive integer, then ${\rm ex}(p+2,B_p)=(p+1)^2/2$ and ${\rm ex}(p+3,B_p)=(p+1)(p+3)/2.$ }

We first characterize the extremal graphs for the two numbers ${\rm ex}(p+2,B_p)$ and ${\rm ex}(p+3,B_p)$ determined in Lemma 3.

{\bf Theorem 4.} {\it Suppose $G$ is a graph of order $p+2$ not containing $B_p.$
\newline\indent (1) If $p$ is even, then $G$ has size ${\rm ex}(p+2,B_p)$ if and only if $G=K_{p+2}-PM$ or $G=\overline{K_1+((p-2)/2)K_2+P_3}.$
\newline\indent (2) If $p$ is odd, then $G$ has size ${\rm ex}(p+2,B_p)$ if and only if $G=K_1\vee (K_{p+1}-PM).$  }

{\bf Proof.} It is easy to check that the three graphs in the theorem satisfy the requirements. Next we prove that they are the only extremal graphs.
Suppose $G$ has size ${\rm ex}(p+2,B_p)$ and its vertices are $v_1,\ldots,v_{p+2}$ with
${\rm deg}_{\overline{G}}(v_1)\le{\rm deg}_{\overline{G}}(v_2)\le\cdots\le {\rm deg}_{\overline{G}}(v_{p+2}).$

(1) $p$ is even. By Lemma 3, $|E(G)|=p(p+2)/2$ which implies that $|E(\overline{G})|=(p+2)/2.$ By the degree-sum formula [4, p.35],
$\sum_{i=1}^{p+2}{\rm deg}_{\overline{G}}(v_i)=2|E(\overline{G})|=p+2,$ implying $\delta(\overline{G})\le 1.$

Case 1. $\delta(\overline{G})=1.$ The degree sequence of $\overline{G}$ is $1,1,\ldots,1.$ Thus $\overline{G}=((p+2)/2)K_2,$ implying that
$G=K_{p+2}-PM.$

Case 2. $\delta(\overline{G})=0.$ By Lemma 1, $v_1$ is the only isolated vertex of $\overline{G}.$ Since $\sum_{i=1}^{p+2}{\rm deg}_{\overline{G}}(v_i)=p+2$
and ${\rm deg}_{\overline{G}}(v_i)\ge 1$ for every $i\ge 2,$ we deduce that $\Delta(\overline{G})=2$ and $\overline{G}$ has a unique vertex of degree $2.$
Thus the degree sequence of $\overline{G}$ is $2,1,\ldots,1,0,$ implying that $\overline{G}=K_1+((p-2)/2)K_2+P_3;$ i.e., $G=\overline{K_1+((p-2)/2)K_2+P_3}.$

(2) $p$ is odd.  By Lemma 3, $|E(G)|=(p+1)^2/2,$ implying $|E(\overline{G})|=(p+1)/2.$ Hence
$\sum_{i=1}^{p+2}{\rm deg}_{\overline{G}}(v_i)=2|E(\overline{G})|=p+1,$
implying $\delta (\overline{G})=0.$ By Lemma 1, $\overline{G}$ has a unique isolated vertex. Thus, the degree sequence of $\overline{G}$ is $1,1,\ldots,1,0.$
It follows that $\overline{G}=K_1+((p+1)/2)K_2;$ i.e, $G=K_1\vee (K_{p+1}-PM).$ $\Box$

{\bf Theorem 5.} {\it Suppose $G$ is a graph of order $p+3$ not containing $B_p.$
\newline\indent (1) If $p=2,$ then $G$ has size ${\rm ex}(p+3,B_p)$ if and only if $G=K_1\vee (2P_2),$ or $G=K_{2,3},$ or $G=\overline{P_5};$
if $p$ is even and $p\ge 4$, then $G$ has size ${\rm ex}(p+3,B_p)$ if and only if $G=(3K_1)\vee (K_p-PM)$ or $G=\overline{P_5}\vee (K_{p-2}-PM).$
\newline\indent (2) If $p$ is odd, then $G$ has size ${\rm ex}(p+3,B_p)$ if and only if $G=K_{p+3}-PM.$  }

{\bf Proof.} It is easy to check that the graphs in the theorem satisfy the requirements. Next we prove that they are the only extremal graphs.
Suppose $G$ has size ${\rm ex}(p+3,B_p)$ and its vertices are $v_1,\ldots,v_{p+3}$ with
${\rm deg}_{\overline{G}}(v_1)\le{\rm deg}_{\overline{G}}(v_2)\le\cdots\le {\rm deg}_{\overline{G}}(v_{p+3}).$

(1) $p$ is even. By Lemma 3,  $|E(G)|=p(p+4)/2$ which implies that $|E(\overline{G})|=(p+6)/2.$ Since
$\sum_{i=1}^{p+3}{\rm deg}_{\overline{G}}(v_i)=2|E(\overline{G})|=p+6,$ we have $\delta(\overline{G})\le 1.$

Case 1. $\delta(\overline{G})=0.$ By Lemma 1, $v_1$ is the only isolated vertex of $\overline{G}$ and hence by Lemma 1 we deduce that
${\rm deg}_{\overline{G}}(v_i)\ge 2$ for every $i\ge 2.$ We have
$$
p+6=\sum_{i=1}^{p+3}{\rm deg}_{\overline{G}}(v_i)\ge 0+(p+2)\times 2=2p+4,
$$
implying $p=2.$ The degree sequence of $\overline{G}$ is $2,2,2,2,0.$ Thus $\overline{G}=K_1+C_4;$ i.e., $G=K_1\vee (2P_2).$

Case 2. $\delta(\overline{G})=1.$ Since $\sum_{i=1}^{p+3}{\rm deg}_{\overline{G}}(v_i)=p+6,$ we have $2\le \Delta (\overline{G})\le 4.$ We assert that
$\Delta (\overline{G})=2.$ If $\Delta (\overline{G})=4$ or $\Delta (\overline{G})=3,$ then the vertex $v_{p+3}$ has at least two neighbors which are leaves in
$\overline{G},$ contradicting the fact that for any $xy\in E(G),$ $|N_{\overline{G}}(x)\cup N_{\overline{G}}(y)|\ge 2.$ Thus the degree sequence of
$\overline{G}$ is $2,2,2,1,1,\ldots,1.$ Note that the two neighbors of a vertex of degree $2$ in $\overline{G}$ cannot be both leaves. It follows that
$\overline{G}=(p/2)K_2+K_3$ or $\overline{G}=((p-2)/2)K_2+P_5;$ i.e., $G=(3K_1)\vee (K_p-PM)$ or $G=\overline{P_5}\vee (K_{p-2}-PM).$
Note that when $p=2,$ these two graphs are $K_{2,3}$ and $\overline{P_5}.$

(2) $p$ is odd. By Lemma 3, $|E(G)|=(p+1)(p+3)/2$ which implies that $|E(\overline{G})|=(p+3)/2.$ Since
$\sum_{i=1}^{p+3}{\rm deg}_{\overline{G}}(v_i)=2|E(\overline{G})|=p+3,$ we have $\delta(\overline{G})\le 1.$ If $\delta(\overline{G})=0,$ then
by Lemma 1 we deduce that for each $i\ge 2,$ ${\rm deg}_{\overline{G}}(v_i)\ge 2,$ which implies
$\sum_{i=1}^{p+3}{\rm deg}_{\overline{G}}(v_i)\ge0+(p+2)\times 2>p+3,$ a contradiction. Hence $\delta(\overline{G})=1,$ which, together with the fact that
$\sum_{i=1}^{p+3}{\rm deg}_{\overline{G}}(v_i)=p+3,$ yields that the degree sequence of $\overline{G}$ is $1,1,\ldots,1.$ It follows that
$\overline{G}=((p+3)/2)K_2;$ i.e., $G=K_{p+3}-PM.$ This completes the proof. $\Box$

In the next two results we determine the Tur\'{a}n number ${\rm ex}(p+4,B_p)$ and characterize the corresponding extremal graphs.

{\bf Theorem 6.} {\it If $p$ is an integer with $p\ge 3,$ then ${\rm ex}(p+4,B_p)=(p+2)(p+3)/2.$ }

{\bf Proof.} Let $G$ be a graph of order $p+4$ and size $e$ with vertices $v_1,\ldots,v_{p+4}$ such that ${\rm deg}(v_i)=d_i,$ $i=1,\ldots,p+4$
and $d_1\ge d_2\ge\cdots\ge d_{p+4.}$ Suppose $e\ge ((p+2)(p+3)/2)+1.$ We will show that $G$ contains $B_p.$ We have
$\sum_{i=1}^{p+4}d_i=2e\ge p^2+5p+8.$ We distinguish two cases.

Case 1. $d_1=p+3.$

If $d_2\le p,$ then $\sum_{i=1}^{p+4}d_i\le p+3+(p+3)p=p^2+4p+3<p^2+5p+8,$ a contradiction. Hence $d_2\ge p+1.$ Then $G$ has $p$ triangles sharing the common
edge $v_1v_2.$

Case 2. $d_1\le p+2.$

If $d_3\le p+1,$ then $\sum_{i=1}^{p+4}d_i\le 2(p+2)+(p+2)(p+1)=p^2+5p+6<p^2+5p+8,$ a contradiction. Hence $d_1=d_2=d_3=p+2.$
Consequently $v_1$ is adjacent to at least one of $v_2$ and $v_3,$ say $v_2.$
Let $v_s\not\in N(v_1)$ and $v_t\not\in N(v_2).$ Then $\{v_3,v_4,\ldots,v_{p+4}\}\setminus \{v_s,v_t\}=N(v_1)\cap N(v_2).$ We have $|N(v_1)\cap N(v_2)|\ge p.$
Thus $G$ has $p$ triangles sharing the common edge $v_1v_2.$ This shows ${\rm ex}(p+4,B_p)\le (p+2)(p+3)/2.$

On the other hand, the graph $G_1=\overline{K_2+C_{p+2}}$ has order $p+4$ and size $(p+2)(p+3)/2,$ and using the criterion in Lemma 1 we check that $G_1$
does not contain $B_p.$ Thus we conclude that ${\rm ex}(p+4,B_p)=(p+2)(p+3)/2.$ $\Box$

Theorem 6 excludes the two small values $p=1,\,2.$ We remark that ${\rm ex}(5,B_1)=6$ agrees with the formula in Theorem 6, but ${\rm ex}(6,B_2)=9$
does not.

{\bf Theorem 7.} {\it Let $p$ be an integer with $p\ge 3$ and let $H$ be a graph of order $p+4$ not containing $B_p.$ Then $H$ has size
${\rm ex}(p+4,B_p)$ if and only if $H=\overline{K_2+C_{i_1}+C_{i_2}+\cdots+C_{i_t}}$ where $i_j\not=4$ for each $j$ and $i_1+i_2+\cdots+i_t=p+2.$ }

{\bf Proof.} It is easy to check that every graph of the form in the theorem satisfies all the requirements; i.e., it is an extremal graph.
We use the criterion in Lemma 1 to check that such a graph does not contain $B_p.$

Conversely suppose that $H$ is a graph of order $p+4$ not containing $B_p$ and $H$ has size ${\rm ex}(p+4,B_p).$ Denote $n=p+4.$
Since $|E(H)|=(p+2)(p+3)/2$ by Theorem 6, $|E(\overline{H})|=n-1.$ By Lemma 1, for any two non-adjacent vertices in $\overline{H},$ we have
$$
|N_{\overline{H}}(x)\cup N_{\overline{H}}(y)|\ge p+4-p-1=3. \eqno (1)
$$

Claim 1. $\overline{H}$ is disconnected.

To the contrary, suppose $\overline{H}$ is connected. Since $\overline{H}$ has order $n$ and size $n-1,$ it is a tree [4, p.68]. Then $\overline{H}$ has two
distinct nonadjacent leaves $u$ and $v,$ since $n\ge 7.$ Now $|N_{\overline{H}}(u)\cup N_{\overline{H}}(v)|=2<3,$ contradicting (1). This proves Claim 1.

Claim 2. Among the components of $\overline{H},$ there is exactly one which is a tree. The tree is $K_2.$

Since $\overline{H}$ has order $n$ and size $n-1,$ it has at least one component which is a tree. We assert that $\overline{H}$ has exactly one such component.
To the contrary, suppose $T_1$ and $T_2$ are two distinct components which are trees. Then $T_1$ has a vertex $s$ and $T_2$ has a vertex $t$ such that
both $s$ and $t$ have degree at most $1.$  We have $|N_{\overline{H}}(s)\cup N_{\overline{H}}(t)|\le 2,$ contradicting (1). Let $T$ be the unique component
of $\overline{H}$ which is a tree. If the order of $T$ is at least $3,$ then as in the above proof of Claim 1 we would obtain a contradiction. Hence $T$
has order at most $2.$ On the other hand, $T$ cannot have order $1.$ Otherwise $T$ is an isolated vertex and $\overline{H}$ has a vertex $w$ other that $T$
which has degree at most $2,$ since the degree sum of $\overline{H}$ is $2(n-1).$ Then we have $|N_{\overline{H}}(T)\cup N_{\overline{H}}(w)|\le 2,$ contradicting (1). It follows that $T=K_2.$

Claim 3. $\overline{H}=K_2+C_{i_1}+C_{i_2}+\cdots+C_{i_t}$ where $i_j\not=4$ for each $j.$

By Claim 2, $\overline{H}=K_2+H_1+H_2+\cdots+H_t$ where $H_j$ is a component of $\overline{H}$ with order $i_j.$ Let $x\in V(K_2)$
and let $y\in V(H_j).$ Then the inequality (1) implies that ${\rm deg}_{\overline{H}}(y)\ge 2.$ Let $v_1,v_2,\ldots, v_n$ be the vertices of
$\overline{H}$ with ${\rm deg}_{\overline{H}}(v_1)={\rm deg}_{\overline{H}}(v_2)=1$ and
$2\le {\rm deg}_{\overline{H}}(v_3)\le\cdots\le {\rm deg}_{\overline{H}}(v_n).$ We have
$2(n-1)=\sum_{i=1}^{n}{\rm deg}_{\overline{H}}(v_i)\ge 2\times 1+(n-2)\times 2=2(n-1),$ implying that ${\rm deg}_{\overline{H}}(v_i)=2$ for each $i=3,\ldots,n.$
Hence $H_j$ is a cycle for every $j=1,2,\ldots,t.$ But none of these cycles can be $C_4.$ To the contrary, suppose $H_s=C_4=abcd.$
Then $|N_{\overline{H}}(a)\cup N_{\overline{H}}(c)|=2,$ contradicting (1). This shows Claim 3 and completes the proof. $\Box$

We determine the Tur\'{a}n number ${\rm ex}(p+5,B_p)$ in the following result.

{\bf Theorem 8.} {\it Let $p$ be a positive integer. Then
$$
{\rm ex}(p+5, B_p)=\begin{cases}\frac{(p+2)(p+5)}{2}\quad {\rm if}\quad p\equiv 1 \,({\rm mod}\, 3)\\
                                \frac{(p+1)(p+6)}{2}\quad {\rm if}\quad p\equiv 0\,{\rm or}\,\, 2\,({\rm mod}\, 3).
                   \end{cases}
$$}

{\bf Proof.} Let $G$ be a graph of order $p+5$ and size $e$ with vertices $v_1,\ldots,v_{p+5}$ such that ${\rm deg}(v_i)=d_i,$ $i=1,\ldots,p+5$
and $d_1\ge d_2\ge\cdots\ge d_{p+5}.$

(1) $p\equiv 1 \,({\rm mod}\, 3).$
Suppose $e\ge ((p+2)(p+5)/2)+1.$ We will show that $G$ contains $B_p.$ We have
$\sum_{i=1}^{p+5}d_i=2e\ge p^2+7p+12.$ We distinguish two cases.

Case 1. $d_1=p+4.$

If $d_2\le p,$ then $\sum_{i=1}^{p+5}d_i\le p+4+(p+4)p=p^2+5p+4<p^2+7p+12,$ a contradiction. Hence $d_2\ge p+1$ and consequently $G$ contains
$p$ triangles sharing the common edge $v_1v_2.$

Case 2. $d_1\le p+3.$

If $d_2\le p+2,$ then $\sum_{i=1}^{p+5}d_i\le p+3+(p+4)(p+2)=p^2+7p+11<p^2+7p+12,$ a contradiction. Thus $d_1=d_2=p+3.$
If $d_3\le p+1,$ then $\sum_{i=1}^{p+5}d_i\le 2(p+3)+(p+3)(p+1)=p^2+6p+9<p^2+7p+12,$ a contradiction. Hence $d_3\ge p+2.$
If $v_1$ and $v_2$ are adjacent, then $G$ contains at least $p+1$ triangles sharing the common edge $v_1v_2.$
If $v_1$ and $v_2$ are nonadjacent, then $G$ contains at least $p$ triangles sharing the common edge $v_1v_3.$ This shows that
${\rm ex}(p+5, B_p)\le (p+2)(p+5)/2.$

On the other hand, the graph $G_2=\overline{((p+5)/3)K_3}$ has order $p+5$ and size $(p+2)(p+5)/2,$ and using Lemma 1 it is easy to verify that $G_2$ does not
contain $B_p.$ Thus, ${\rm ex}(p+5, B_p)=(p+2)(p+5)/2.$

(2) $p\equiv 0\,{\rm or}\,\, 2\,({\rm mod}\, 3).$ Suppose $e\ge ((p+1)(p+6)/2)+1.$ We will show that $G$ contains $B_p.$ We have
$\sum_{i=1}^{p+5}d_i=2e\ge p^2+7p+8.$ We distinguish three cases.

Case 1. $d_1=p+4.$

The proof for this case is similar to Case 1 of (1) above.

Case 2. $d_1=p+3.$

If $d_3\le p+1,$ then $\sum_{i=1}^{p+5}d_i\le 2(p+3)+(p+3)(p+1)=p^2+6p+9<p^2+7p+8,$ a contradiction. Hence $d_3\ge p+2.$ There are two possible
values of $d_2.$ If $d_2=p+3,$ as in Case 2 of (1) above we can show that $G$ contains $B_p.$ If $d_2=p+2,$ we have $d_2=d_3=p+2.$ Since
$d_1=p+3,$ $v_1$ is adjacent to at least one of $v_2$ and $v_3.$ Without loss of generality, suppose  $v_1$ is adjacent to $v_2.$ Now $G$ contains
$p$ triangles sharing the common edge $v_1v_2.$

Case 3. $d_1\le p+2.$

If $d_{p+3}\le p+1,$ then $\sum_{i=1}^{p+5}d_i\le (p+2)(p+2)+3(p+1)=p^2+7p+7<p^2+7p+8,$ a contradiction. Hence $d_1=d_2=\cdots=d_{p+3}=p+2.$

Subcase 3.1. $d_{p+5}=p+2.$ Now ${\rm deg}_{\overline{G}}(v_1)=\cdots={\rm deg}_{\overline{G}}(v_{p+5})=2.$ Thus $\overline{G}$ is a cycle
or a union of vertex-disjoint cycles. The condition $p\equiv 0\,{\rm or}\,\, 2\, ({\rm mod}\, 3)$ implies that the order $p+5$ is not
divisible by $3.$ It follows that $\overline{G}$ has a component which is a cycle of length at least $4.$  By Lemma 2, $G$ contains $B_p.$

Subcase 3.2. $d_{p+5}\le p+1.$ The condition $\sum_{i=1}^{p+5}d_i\ge p^2+7p+8$ yields that $d_{p+4}\ge p+1$ and that if $d_{p+4}=p+1$ then $d_{p+5}=p+1.$
First suppose in $G,$ every non-neighbor of $v_{p+5}$ has degree $p+2.$ Since ${\rm deg}_{\overline{G}}(v_{p+5})\ge 3,$ we deduce that $v_{p+5}$
has two adjacent non-neighbors, say $v_1$ and $v_2,$ in $G.$ Then $|N_{\overline{G}}(v_1)\cup N_{\overline{G}}(v_2)|\le 3<4.$ By Lemma 1, $G$ contains $B_p.$

It remains to consider the case when in $G,$ $v_{p+5}$ has at least one non-neighbor whose degree is not $p+2.$
Such a non-neighbor can only be $v_{p+4},$ since $d_1=d_2=\cdots=d_{p+3}=p+2.$ Then the degree sum condition implies that
$d_{p+4}=d_{p+5}=p+1.$ It is easy to deduce that $\overline{G}$ contains a component $F$ of order $6$ depicted in Figure 2
where the two cut-vertices are $v_{p+4}$ and $v_{p+5}.$
\vskip 3mm
\par
 \centerline{\includegraphics[width=1.8in]{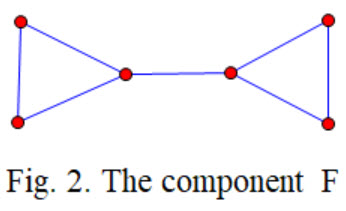}}
\par
The graph $\overline{G}-V(F)$ is a regular graph of degree $2$ and has order $p-1$ which is not divisible by $3.$ It follows that
$\overline{G}-V(F)$ and hence $\overline{G}$ contains a component which is a cycle of length at least $4.$ By Lemma 2 we deduce that $G$
contains $B_p.$

So far we have proved that ${\rm ex}(p+5, B_p)\le (p+1)(p+6)/2.$ Next we show that this upper bound can be attained.

First suppose $p\equiv 0 \,({\rm mod}\, 3).$ Let $W$ denote the graph of order $8$ in Figure 3.
\vskip 3mm
\par
 \centerline{\includegraphics[width=2.1in]{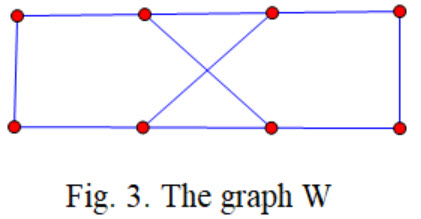}}
\par
Denote $G_3=\overline{((p-3)/3)K_3+W}.$ Then $G_3$ has order $p+5$ and size $(p+1)(p+6)/2,$ and by using Lemma 1 we can verify that $G_3$ does not contain
$B_p.$

Next suppose $p\equiv 2 \,({\rm mod}\, 3).$ Denote $G_4=\overline{((p+1)/3)K_3+K_4}.$ Then $G_4$ has order $p+5$ and size $(p+1)(p+6)/2,$ and $G_4$ does not contain $B_p.$ This completes the proof. $\Box$

We will need the following three lemmas to characterize the extremal graphs for the Tur\'{a}n number ${\rm ex}(p+5, B_p).$

{\bf Lemma 9.} {\it Let $G$ be a graph of order $p+5$ not containing $B_p.$ If $x$ and $y$ are two distinct vertices with
${\rm deg}_{\overline{G}}(x)={\rm deg}_{\overline{G}}(y)=2$ and the distance between them in $\overline{G}$ is at most $2,$ then $x$ and $y$
are adjacent in $\overline{G}.$ }

{\bf Proof.} To the contrary, suppose that $x$ and $y$ are nonadjacent in $\overline{G}.$ Then by Lemma 1,
$|N_{\overline{G}}(x)\cup N_{\overline{G}}(y)|\ge(p+5)-p-1=4.$ Also, the condition that the distance between $x$ and $y$ in $\overline{G}$ is at most $2$
implies that $N_{\overline{G}}(x)\cap N_{\overline{G}}(y)\not=\phi.$ We have
$$
|N_{\overline{G}}(x)\cup N_{\overline{G}}(y)|=|N_{\overline{G}}(x)|+|N_{\overline{G}}(y)|-|N_{\overline{G}}(x)\cap N_{\overline{G}}(y)|\le 2+2-1=3,
$$
a contradiction. $\Box$

{\bf Lemma 10.} {\it Let $G$ be a graph of order $p+5$ not containing $B_p.$ If $P$ is a path in $\overline{G}$ and every internal vertex of $P$
has degree $2$ in $\overline{G},$ then the length of $P$ is at most $3.$ }

{\bf Proof.} To the contrary, suppose that $P=v_1v_2\ldots v_k$ where $k\ge 5$ and in $\overline{G},$ each $v_i$ has degree $2$ for $i=2,3,\ldots,k-1.$
But then
$$
|N_{\overline{G}}(v_2)\cup N_{\overline{G}}(v_4)|=3<4=(p+5)-p-1,
$$
contradicting Lemma 1. $\Box$

{\bf Lemma 11.} {\it Let $G$ be a graph of order $p+5$ and size ${\rm ex}(p+5, B_p)$ not containing $B_p$ where $p\ge 4.$
\newline\indent (1) If $p\equiv 1 \,({\rm mod}\, 3),$ then $\overline{G}$ is $2$-regular;
\newline\indent (2) if $p\equiv 0\,{\rm or}\,\, 2\, ({\rm mod}\, 3),$ then the degree sequence of $\overline{G}$ is $3,3,3,3,2,\ldots,2.$ }

{\bf Proof.} Let $v_1,\ldots,v_{p+5}$ be the vertices of $G$ with ${\rm deg}_{\overline{G}}(v_i)=d_i,$ $i=1,\ldots,p+5$ such that
 $d_1\ge d_2\ge\cdots\ge d_{p+5}.$

 (1) $p\equiv 1 \,({\rm mod}\, 3).$ By Theorem 8, $|E(G)|=(p+2)(p+5)/2,$ implying that $\overline{G}$ has size $p+5.$ If $\overline{G}$ has
 an isolated vertex $x,$ then by Lemma 1 we deduce that for any vertex $u\in V(G)\setminus\{x\},$ ${\rm deg}_{\overline{G}}(u)\ge 4.$ Hence
 $2p+10=2|E(\overline{G})|=\sum_{i=1}^{p+5}{\rm deg}_{\overline{G}}(v_i)\ge 0+(p+4)\times 4=4p+16,$ a contradiction. If $\overline{G}$ has a leaf
 $y,$ let $v$ be the neighbor of $y$ in $\overline{G}.$ Then by Lemma 1, for any vertex $w\in V(G)\setminus\{y,v\},$ ${\rm deg}_{\overline{G}}(w)\ge 3.$
 We have $2p+10=2|E(\overline{G})|=\sum_{i=1}^{p+5}{\rm deg}_{\overline{G}}(v_i)\ge 1+1+(p+3)\times 3=3p+11,$ a contradiction. Thus
 $\delta(\overline{G})\ge 2.$ Now the inequality $2p+10=2|E(\overline{G})|=\sum_{i=1}^{p+5}{\rm deg}_{\overline{G}}(v_i)\ge (p+5)\times 2=2p+10$
 forces ${\rm deg}_{\overline{G}}(v_i)=2$ for each $i=1,\ldots, p+5;$ i.e., $\overline{G}$ is $2$-regular.

 (2) $p\equiv 0\,{\rm or}\,\, 2\, ({\rm mod}\, 3).$ By Theorem 8, $\overline{G}$ has size $p+7.$ As in the above proof of part (1), we deduce that
 $\delta(\overline{G})\ge 2.$ If $\delta(\overline{G})\ge 3,$ then $2p+14=2|E(\overline{G})|=\sum_{i=1}^{p+5}{\rm deg}_{\overline{G}}(v_i)\ge (p+5)\times 3 =3p+15,$ a contradiction. Hence $\delta(\overline{G})=2.$ From $2p+14=2|E(\overline{G})|=\sum_{i=1}^{p+5}{\rm deg}_{\overline{G}}(v_i)\ge (p+4)\times 2+
 \Delta(\overline{G})$ we obtain $\Delta(\overline{G})\le 6.$

If $\Delta(\overline{G})=6,$ then the degree sequence of $\overline{G}$ is $6,2,\ldots,2;$ if $\Delta(\overline{G})=5,$ then the degree sequence of $\overline{G}$ is $5,3,2,\ldots,2;$ if $\Delta(\overline{G})=4,$ then the degree sequence of $\overline{G}$ is $4,4,2,\ldots,2$ or $4,3,3,2,\ldots,2.$
In the first three cases, i.e., $\overline{G}$ has degree sequence $6,2,\ldots,2,$ or $5,3,2,\ldots,2,$ or $4,4,2,\ldots,2,$ $v_1$ has at least three
neighbors of degree $2$ in $\overline{G}$ which are pairwise adjacent by Lemma 9, a contradiction. Now suppose $\overline{G}$ has degree sequence
$4,3,3,2,\ldots,2.$ We will show that this cannot occur.

By the above argument, $v_1$ has at most two neighbors of degree $2$ in $\overline{G}.$ Hence $v_2,v_3\in N_{\overline{G}}(v_1).$ Without loss of generality,
suppose $N_{\overline{G}}(v_1)=\{v_2,v_3,v_4,v_5\}.$ By Lemma 9, $v_4$ and $v_5$ are adjacent in $\overline{G}.$ We distinguish two cases.

Case 1. $v_2$ and $v_3$ are adjacent in $\overline{G}.$

Subcase 1.1. $|N_{\overline{G}}(v_2)\cap N_{\overline{G}}(v_3)|\ge 2.$ Since $d_2=d_3=3,$ we have $|N_{\overline{G}}(v_2)\cap N_{\overline{G}}(v_3)|=2.$
Without loss of generality, suppose $N_{\overline{G}}(v_2)\cap N_{\overline{G}}(v_3)=\{v_1,v_6\}.$ Then $H=\overline{G}[v_1,v_2,v_3,v_4,v_5,v_6]$ is a component
of $\overline{G}$ and $\overline{G}-V(H)$ is a $2$-regular graph of order $p-1.$ Since $3$ does not divide $p-1,$ $\overline{G}$ contains a component
which is a cycle of length at least $4.$ By Lemma 2, $G$ contains $B_p,$ a contradiction.

Subcase 1.2. $|N_{\overline{G}}(v_2)\cap N_{\overline{G}}(v_3)|=1.$ Without loss of generality, suppose $N_{\overline{G}}(v_2)=\{v_1,v_3,v_6\}$ and
$N_{\overline{G}}(v_3)=\{v_1,v_2,v_7\}.$ By Lemma 10, $v_6$ and $v_7$ are adjacent in $\overline{G}.$ Now  $v_3$ and $v_6$ are adjacent in $G,$
and $|N_{\overline{G}}(v_3)\cup N_{\overline{G}}(v_6)|=3<4=(p+5)-p-1.$ By Lemma 1, $G$ contains $B_p,$ a contradiction.

Case 2.  $v_2$ and $v_3$ are nonadjacent in $\overline{G}.$

Without loss of generality, suppose $N_{\overline{G}}(v_2)=\{v_1,v_6,v_7\}$ and $N_{\overline{G}}(v_3)=\{v_1,v_8,v_9\}.$ By Lemma 9, $v_6$ and $v_7$
are adjacent in $\overline{G},$ and $v_8$ and $v_9$ are adjacent in $\overline{G}.$ Thus, $R=\overline{G}[v_1,v_2,\ldots,v_9]$ is a component of
$\overline{G},$ and $\overline{G}-V(R)$ is a $2$-regular graph of order $p-4.$ Since $p-4$ is not divisible by $3,$ $\overline{G}$ has a component
which is a cycle of length at least $4.$ By Lemma 2, $G$ contains $B_p,$ a contradiction.

We have proved that $\Delta(\overline{G})\le 3.$ Since $\delta(\overline{G})=2$ and $\overline{G}$ has size $p+7,$ we deduce that the
degree sequence of $\overline{G}$ is $3,3,3,3,2,\ldots,2.$ $\Box$

We denote by $Q$ the graph of order $10$ in Figure 4.
\vskip 3mm
\par
 \centerline{\includegraphics[width=1.8in]{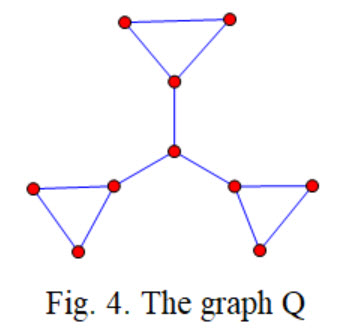}}
\par
We characterize the extremal graphs for the Tur\'{a}n number ${\rm ex}(p+5, B_p)$ in the following result.

{\bf Theorem 12.} {\it Let $p$ be an integer with $p\ge 4$ and let $G$ be a graph of order $p+5$ not containing $B_p.$ Then $G$ has size
${\rm ex}(p+5, B_p)$ if and only if
\newline\indent (1) $G=\overline{((p-3)/3)K_3+W}$ where $W$ is the graph depicted in Figure 3 if $p\equiv 0 \,({\rm mod}\, 3);$
\newline\indent (2) $G=\overline{((p+5)/3)K_3}$ if $p\equiv 1 \,({\rm mod}\, 3);$
\newline\indent (3) $G=\overline{((p+1)/3)K_3+K_4}$ or $G=\overline{((p-5)/3)K_3+Q}$ where $Q$ is the graph depicted in Figure 4
if $p\equiv 2 \,({\rm mod}\, 3).$ }

{\bf Proof.} The sufficiency can be easily verified by using Lemma 1. Now we prove the necessity. Suppose $G$ has size ${\rm ex}(p+5, B_p).$
We will consider the complement graph $\overline{G}.$

We first prove part (2) where $p\equiv 1 \,({\rm mod}\, 3).$ By Lemma 11, $\overline{G}$ is $2$-regular, or equivalently $\overline{G}$ is a union
of vertex-disjoint cycles.  By Lemma 2, every cycle must be a triangle. Thus $\overline{G}=((p+5)/3)K_3.$

Next we treat parts (1) and (3). Suppose $p\equiv 0\,\, {\rm or}\,\, 2 \,({\rm mod}\, 3).$ By Lemma 11, the degree sequence of $\overline{G}$
is $3,3,3,3,2,\ldots,2.$ By Lemma 9, in $\overline{G}$ each vertex of degree $3$ must have at least one neighbor of degree $3.$

Claim. The four vertices of degree $3$ lie in one component of $\overline{G}$ and there is a vertex of degree $3$ with at least two neighbors of
degree $3.$

To the contrary, suppose the four vertices of degree $3$ do not lie in one component of $\overline{G}.$ Then they lie in two components and by Lemma 9
we deduce that each of these two components is the graph $F$ of order $6$ in Figure 2. Every component of $\overline{G}$ other than these two components
is a cycle and at least one of these cycles has length at least $4,$ since the order $p+5$ of $\overline{G}$ is not divisible by $3.$ But this contradicts
Lemma 2. Hence the first conclusion in the Claim is proved. The second one can be proved similarly.

Denote by $Z$ the component of $\overline{G}$ in which the four vertices of degree $3$ lie. By the above properties,  $Z$ contains either a claw or a path of order $4$ each of whose vertices has degree $3.$ In both cases, by using Lemmas 9, 10 and 2 it can be verified that $Z$ has order at most $10.$ By Lemma 2, the graph $\overline{G}-V(Z)$ is a union of vertex-disjoint triangles. Hence $p+5=|\overline{G}|\equiv |Z|\,({\rm mod}\, 3).$

If $p\equiv 0\,({\rm mod}\, 3),$ $|Z|\in \{5,8\};$ if $p\equiv 2\,({\rm mod}\, 3),$ $|Z|\in \{4,7,10\}.$

Suppose $p\equiv 0\,({\rm mod}\, 3).$ By Lemma 1, $|Z|\not=5.$ If $|Z|=8,$ then $Z$ is the graph $W$ in Figure 3. It follows that
$\overline{G}=((p-3)/3)K_3+W.$ This proves part (1).

Suppose $p\equiv 2\,({\rm mod}\, 3).$ By Lemma 1, $|Z|\not= 7.$ If $|Z|=4,$ then $Z=K_4$  and if $|Z|=10,$ then $Z=Q.$ This proves part (3). The proof is complete. $\Box$

Theorem 12 does not include the small values $p=1,2,3.$ The information about them is as follows.
$K_{3,3}$ is the unique extremal graph for ${\rm ex}(6, B_1)=9.$ $K_{3,4}$ is the unique
extremal graph for ${\rm ex}(7, B_2)=12.$ There are two extremal graphs for ${\rm ex}(8, B_3)=18$ which are depicted in Figure 5.
\vskip 3mm
\par
 \centerline{\includegraphics[width=3.1in]{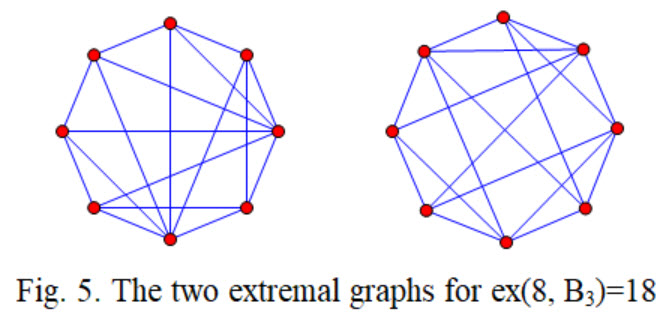}}
\par

We determine the Tur\'{a}n number ${\rm ex}(p+6,B_p)$ in the following result.

{\bf Theorem 13.} {\it Let $p$ be a positive integer.
\newline\indent (1) If $p$ is odd and $p\ge 5,$ then ${\rm ex}(p+6,B_p)=(p+3)(p+5)/2;$
\newline\indent (2) if $p$ is even and $p\not=2,6,10,$ then ${\rm ex}(p+6,B_p)=1+(p+2)(p+6)/2.$ }

{\bf Proof.} Let $G$ be a graph of order $p+6$ and size $e$ with vertices $v_1,\ldots,v_{p+6}$ such that ${\rm deg}(v_i)=d_i,$ $i=1,\ldots,p+6$
and $d_1\ge d_2\ge\cdots\ge d_{p+6}.$

(1) $p$ is odd and $p\ge 5.$ Suppose $e\ge ((p+3)(p+5)/2)+1.$ We will show that $G$ contains $B_p.$ We have
$\sum_{i=1}^{p+6}d_i=2e\ge p^2+8p+17.$ We distinguish three cases.

Case 1. $d_1=p+5.$

If $d_2\le p,$ then $\sum_{i=1}^{p+6}d_i\le p+5+(p+5)\times p=p^2+6p+5<p^2+8p+17,$ a contradiction. Hence $d_2\ge p+1$ and $G$ contains $p$ triangles
sharing the common edge $v_1v_2.$

Case 2. $d_1=p+4.$

If $d_3\le p+2,$ then $\sum_{i=1}^{p+6}d_i\le 2(p+4)+(p+4)(p+2)=p^2+8p+16<p^2+8p+17,$ a contradiction. Hence $d_3\ge p+3.$ $v_1$ has at least one neighbor
in $\{v_2,v_3\}.$ If $v_1$ is adjacent to $v_2,$ then $|N(v_1)\cap N(v_2)|\ge p+1$ and thus $G$ contains $p+1$ triangles sharing the common edge
$v_1v_2;$ if $v_1$ is adjacent to $v_3,$ then similarly we see that $G$ contains $p$ triangles sharing the common edge $v_1v_3.$

Case 3. $d_1\le p+3.$

If $d_4\le p+2,$ then $\sum_{i=1}^{p+6}d_i\le 3(p+3)+(p+3)(p+2)=p^2+8p+15<p^2+8p+17,$ a contradiction. Hence $d_4\ge p+3$ and $d_1=d_2=d_3=d_4=p+3.$
$v_1$ has at least one neighbor in $\{v_2,v_3,v_4\}.$ Without loss of generality, suppose $v_1$ and $v_2$ are adjacent. Then
$|N(v_1)\cap N(v_2)|\ge p$ and $G$ contains $p$ triangles sharing the common edge $v_1v_2.$

We have proved that ${\rm ex}(p+6,B_p)\le (p+3)(p+5)/2.$ Next we construct graphs to show that this upper bound can be attained.
First suppose $p\equiv 1\,(\rm mod)\, 4.$ The graph $G_5=\overline{K_3+((p+3)/4)K_4}$ has order $p+6,$ size $(p+3)(p+5)/2$ and by using Lemma 1
we can verify that $G_5$ does not contain $B_p.$ Now suppose $p\equiv 3\,(\rm mod)\, 4.$ Denote by $Y$ the graph of order $10$ in Figure 6.
\vskip 3mm
\par
 \centerline{\includegraphics[width=1.5in]{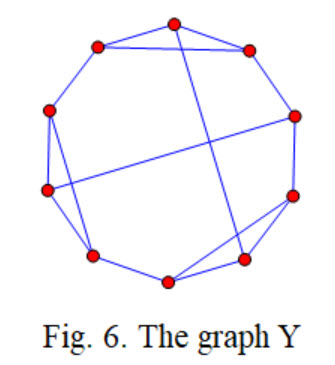}}
\par
\noindent Then the graph $G_6=\overline{K_3+((p-7)/4)K_4+Y}$ has order $p+6$ and size $(p+3)(p+5)/2,$ and $G_6$ does not contain $B_p.$ This proves
${\rm ex}(p+6,B_p)=(p+3)(p+5)/2.$

(2) $p$ is even and $p\not=2,6,10.$ Suppose $e\ge 1+((p+2)(p+6)/2)+1.$ We will show that $G$ contains $B_p.$ We have
$\sum_{i=1}^{p+6}d_i=2e\ge p^2+8p+16.$ We distinguish three cases.

Case 1. $d_1=p+5.$

If $d_2\le p,$ then $\sum_{i=1}^{p+6}d_i\le p+5+(p+5)p=p^2+6p+5<p^2+8p+16,$ a contradiction. Thus $d_2\ge p+1$ and $G$ contains $p$ triangles sharing
the common edge $v_1v_2.$

Case 2. $d_1=p+4.$

If $d_3\le p+1,$ then $\sum_{i=1}^{p+6}d_i\le 2(p+4)+(p+4)(p+1)=p^2+7p+12<p^2+8p+16,$ a contradiction. Hence $d_2\ge d_3\ge p+2.$
$v_1$ has at least one neighbor in $\{v_2,v_3\}$ and consequently $G$ contains $p$ triangles sharing the common edge $v_1v_2$ or $v_1v_3.$

Case 3. $d_1\le p+3.$

If $d_4\le p+2,$ then $\sum_{i=1}^{p+6}d_i\le 3(p+3)+(p+3)(p+2)=p^2+8p+15<p^2+8p+16,$ a contradiction. Hence $d_4\ge p+3.$ In fact,
$d_1=d_2=d_3=d_4=p+3.$ It follows that $G$ contains $B_p.$

We have proved that ${\rm ex}(p+6,B_p)\le 1+(p+2)(p+6)/2.$ Next we show that this upper bound can be attained.

First suppose $p\equiv 0\,({\rm mod}\, 4).$ Denote by $S$ the graph of order $7$ in Figure 7.
\vskip 3mm
\par
 \centerline{\includegraphics[width=1.6in]{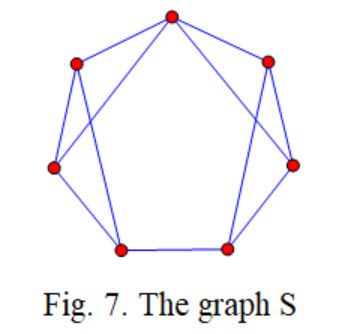}}
\par
\noindent Then the graph $G_7=\overline{K_3+((p-4)/4)K_4+S}$ has order $p+6$ and size $1+(p+2)(p+6)/2,$ and by using Lemma 1 we can verify that $G_7$ does not contain $B_p.$

Now suppose $p\equiv 2\,({\rm mod}\, 4).$ We have $p\ge 14.$ The graph
$$
G_8=\overline{K_3+((p-14)/4)K_4+S+Y}
$$
has order $p+6$ and size $1+(p+2)(p+6)/2,$ and $G_8$ does not contain $B_p.$ This completes the proof. $\Box$

Theorem 13 does not include the small values $p=2,3,6,10.$ The information about them is as follows. ${\rm ex}(8,B_2)=16,$ ${\rm ex}(9,B_3)=21,$
${\rm ex}(12,B_6)=48$ and ${\rm ex}(16,B_{10})=96.$

It seems difficult to characterize the extremal graphs for ${\rm ex}(p+6,B_p).$ A computer search shows that there are $16$ extremal graphs
for ${\rm ex}(9,B_3)=21$ and there are $20$ extremal graphs for ${\rm ex}(12,B_6)=48.$

The results in this paper suggest that perhaps there is no uniform formula for the  Tur\'{a}n number ${\rm ex}(n, B_p)$ for all pairs $(n,p).$

\vskip 5mm
{\bf Acknowledgement.} This research  was supported by the NSFC grants 11671148 and 11771148 and Science and Technology Commission of Shanghai Municipality (STCSM) grant 18dz2271000.

\end{document}